\documentclass{amsart}
\usepackage{mathrsfs}

\newtheorem{Theorem}{Theorem}
\newtheorem{Lemma}{Lemma}
\newtheorem{Corollary}{Corollary}

\def\C{\mathbb C}
\def\H{\mathbb H}
\def\P{\mathbb P}
\def\R{\mathbb R}
\def\Z{\mathbb Z}
\def\CC{\mathscr C}
\def\Re{\operatorname{Re}}
\def\Mod{\,\operatorname{mod}\,}

\begin{document}

\title{Ap\'ery limits and special values of $L$-functions}
\author{Yifan Yang}
\address{Department of Applied Mathematics \\
  National Chiao Tung University \\
  Hsinchu 300, Taiwan}
\email{yfyang@math.nctu.edu.tw}
\date\today
\begin{abstract} We describe a general method to determine the
  Ap\'ery limits of a differential equation that have a
  modular-function origin. As a by-product of our analysis, we
  discover a family of identities involving the special values of
  $L$-functions associated with modular forms. The proof of these
  identities is independent of differential equations and Ap\'ery
  limits.
\end{abstract}


\maketitle

\begin{section}{Introduction} \label{section: introduction} In 1978,
  R.~Ap\'ery proved that
  $\zeta(3)=\sum_{n=1}^\infty n^{-3}$ is an irrational number by
  constructing two sequences
\begin{equation*}
  a_n=\sum_{k=0}^n\binom nk^2\binom{n+k}k^2,
\end{equation*}
$$
  b_n=\sum_{k=0}^n\binom nk^2\binom{n+k}k^2\left\{
    \sum_{m=1}^n\frac1{m^3}+\sum_{m=1}^k\frac{(-1)^{m-1}}
    {\displaystyle 2m^3\binom nm\binom{n+m}m}\right\},
$$
and then showing that $b_n/a_n$ converges to $\zeta(3)$ fast enough to
ensure irrationality of $\zeta(3)$ (see \cite{Poorten}).
Another remarkable discovery of
Ap\'ery is that $a_n$ and $b_n$ satisfy the recursive relation
$$
  (n+2)^3u_{n+2}-(34n^3+153n^2+231n+117)u_{n+1}+(n+1)^3 u_n=0 \
  (u_n=a_n\text{ or }b_n).
$$
Thus, if we set $A(t)=\sum_{n=0}^\infty a_nt^n$ and
$B(t)=\sum_{n=0}^\infty b_nt^n$, then the functions $A(t)$ and $B(t)$
satisfy the differential equations
\begin{equation}
\label{Apery's differential equation}
  (1-34t+t^2)\theta^3A+(3t^2-51t)\theta^2A+(3t^2-27t)\theta A+(t^2-5t)A=0
\end{equation}
and
$$
  (1-34t+t^2)\theta^3B+(3t^2-51t)\theta^2B+(3t^2-27t)\theta B+(t^2-5t)B=6t,
$$
where $\theta$ denotes the differential operator $td/dt$.

Ap\'ery also had an analogous result for $\zeta(2)=\pi^2/6$. He
showed that if $\{a_n\}$ and $\{b_n\}$ are sequences of rational
numbers satisfying the recursive relation
$$
  (n+2)^2u_{n+2}-(11n^2+33n+25)u_{n+1}-(n+1)^2u_n=0 \quad
  (u_n=a_n\text{ or }b_n)
$$
with the initial values
$$
  a_{-1}=0, \quad a_0=1, \quad b_0=0, \quad b_1=5,
$$
then $b_n/a_n$ converges to $\zeta(2)$. Again, due to the recursive
relation, the generating functions $A(t)=\sum a_nt^n$ and $B(t)=\sum
b_nt^n$ satisfy
\begin{equation} \label{Apery's differential equation 2}
  (1-11t-t^2)\theta^2A-(11t+2t^2)\theta A-(3t+t^2)A=0
\end{equation}
and
$$
  (1-11t-t^2)\theta^2B-(11t+2t^2)\theta B-(3t+t^2)B=5t,
$$
respectively.

A common feature of Ap\'ery's two examples is the existence of a
differential equation $LA(t)=0$ with regular singularities whose local
exponents at $t=0$ are all $0$ such that if $A(t)=1+\sum a_nt^n$ is
the unique holomorphic solution at $t=0$ and $B(t)=t+\sum b_nt^n$ is
the unique holomorphic solution of the inhomogeneous differential
equation $LB(t)=t$ at $t=0$, then the ratios $b_n/a_n$ converge to a
special value of the Riemann zeta function. Inspired by these two
examples, Zudilin et al. \cite{ASZ,Zudilin1,Zudilin3} considered the
Ap\'ery limits of a differential equation. The general setting is
described as follows.

Let $Lf(t)=0$ be a linear differential equation with polynomial
coefficients and regular singularities. Assume that the local
exponents of $L$ at $t=0$ are all $0$ so that the monodromy around
$t=0$ is maximally unipotent. In other words, the differential
operator $L$ takes the form
$$
  \theta^k+tP_1(\theta)+t^2P_2(\theta)+\cdots+t^dP_d(\theta),
  \qquad \theta=td/dt,
$$
where $P_j$ are polynomials of degree $\le k$. Let $A(t)=1+\sum
a_nt^n$ be the unique holomorphic solution at $t=0$. Then the sequence
$\{a_n\}$ satisfy a $(d+1)$-term recursive relation
$$
  (n+d)^ka_{n+d}+P_1(n+d-1)a_{n+d-1}+\cdots+P_d(n)a_n=0
$$
with initial values $a_{-d+1}=\cdots=a_{-1}=0$ and $a_0=1$. Now for
each integer $j$ from $1$ to $d-1$, we let
$B(t)=t^j+\sum_{n=j+1}^\infty b_nt^n$ be a solution of the
inhomogeneous differential equation $LB(t)=j^kt^j$. (Note that when
$j\ge d$, there may not exist a solution of $LB(t)=t^j$ that is
holomorphic at $t=0$.) The coefficients
$b_n$ also satisfy the recursive relation
$$
  (n+d)^kb_{n+d}+P_1(n+d-1)b_{n+d-1}+\cdots+P_d(n)b_n=0
$$
with initial values $b_{j-d+1}=\cdots=b_{j-1}=0$ and $b_j=1$. Then the
$j$th \emph{Ap\'ery limit} is defined to be the limit of $b_n/a_n$.
For example, in \cite{Zudilin1}, Zudilin gave a sixth order
differential equation whose Ap\'ery limit is the Catalan number
$\sum_{n=0}^\infty(-1)^n/(2n+1)^2$ and a fifth order differential
equation whose Ap\'ery limit is $\zeta(4)$. In \cite{Zudilin3}, he
also found a third order differential equations whose first and second
Ap\'ery limits give simultaneous approximations to $\log 2$ and
$\zeta(2)$. The paper also contains an example of a sixth order
differential equation that gives simultaneous approximation to
$\zeta(2)$ and $\zeta(3)$. More recently, Almkvist et al. \cite{ASZ}
consider the Ap\'ery limits of fourth order differential equations of
Calabi-Yau type. The numerical computation finds that the recognized
limits are all rational combinations of values of zeta functions or
Dirichlet series associated with odd characters modulo $3$ or $4$.

In practice, one would be interested only in the cases where the
coefficients $a_n$ are integers having nice properties, such as having
a closed form in terms of binomial coefficients. When such differential
equations have order $2$ or $3$, they are often related to modular
forms. This is because in these cases the monodromy groups happen to
be nice arithmetic subgroups of $SL(2,\R)$. Therefore, it is possible
to use the theory of modular forms and modular functions to determine
the Ap\'ery limits for these cases. In fact, this has been done by
Beukers \cite{Beukers} earlier. However, we remark that in
\cite{Beukers} the main goal is to give a modular-function
interpretation of Ap\'ery's irrationality results, so the argument is
not readily applicable to general situations. Thus, the main purpose
of the present article is to present a method that works for more
general differential equations. In Section \ref{section: approach} we
will briefly review Beukers' argument, and then describe our general
approach. In Sections \ref{section: elliptic} and \ref{section: cusp},
we specialize our method to the cases where the non-zero singularity
closest to the origin corresponds to an elliptic point or a cusp of
the underlying modular curve. The examples we work out include cases
($e$), ($h$), and ($\beta$) of \cite{ASZ}. These are some of the
examples where \cite{ASZ} fails to determine their Ap\'ery limits. We
also give two examples in which the Ap\'ery limits are values of
$L$-functions associated with cusp forms.

As a by-product, our analysis in one particular example leads us to
a family of identities involving the special values of $L$-functions
associated with modular forms of weight $3$. (See Lemma \ref{lemma: L
  values} in Section \ref{section: L values}.) A typical example is
$$
  \sum_{n\equiv 1\Mod 12}\frac{c_n}{n^2}=\frac{2+\sqrt3}3
  \sum_{n=1}^\infty\frac{c_n}{n^2}
$$
for the cusp form
$f(\tau)=\eta(2\tau)^3\eta(6\tau)^3=\sum_{n=1}^\infty c_ne^{2\pi
in\tau}$ of weight $3$. The proof of these identities is independent of
differential equations and Ap\'ery limits. We expect that our argument
can be extended to modular forms of higher weight, but we do not
attempt to do so since this will be too far astray from the main topic
of the paper.

After finishing this paper, Zudilin kindly informed of us that in an
unpublished manuscript \cite{Zagier}, Zagier also has described an
approach to the determination of Ap\'ery limits very close to ours
(presumably in the same spirit as \cite{Beukers}), but we do not
know to what extent the two papers overlap.
\medskip

\centerline{\sc Acknowledgments}
\medskip

The author would like to thank Professor Hildebrand of the University
of Illinois for providing references to Tauberian theorems. He would
also like to thank Professors Yui and Zudilin for their interest in
the work and providing valuable comments. The work
was done while the author was a visiting scholar at the Queen's
University, Canada. The author would like to thank Professor Yui and
the Queen's University for their warm hospitality. The visit was
supported by Grant 96-2918-I-009-005 of the National Science Council
of Taiwan and Discovery Grant of Professor Yui of the National
Sciences and Engineering Council of Canada.
\end{section}

\begin{section}{Modular-function approach} \label{section: approach}
Recall that a result of Stiller \cite{Stiller} (see also \cite{Yang})
states that if $A(\tau)$ is a (meromorphic) modular form of weight $k$
with character $\chi$ and $t(\tau)$ is a non-constant modular function
on a subgroup $\Gamma$ of $SL(2,\R)$ commensurable with $SL(2,\Z)$,
then $A,~\tau A,\dots,~\tau^k A$, as functions of $t$, are
linearly independent solutions of a $(k+1)$-st order linear
differential equation $LA=0$ with algebraic functions of $t$ as
coefficients. We assume that the differential equation $LA=0$ has
polynomial coefficients. This is true whenever $t(\tau)$ is a
uniformizer (Hauptmodul) of the modular curve $X(\Gamma)$. From now on
we assume that $\Gamma$ is of genus zero so that such a uniformizer
exists.

We assume that $t(i\infty)=0$ and the Fourier expansion of $A(\tau)$
starts from $1$. Then the differential equation satisfied by $A$ and
$t$ has local exponents all equal to $0$ at $t=0$. When the function
$t:\tau\mapsto t(\tau)$ is locally one-to-one and $A(\tau)$ is
holomorphic at a point $\tau_0$, the function $A(\tau)$ is a
well-defined analytic function of $t$ near $t_0=t(\tau_0)$.
Conversely, if either $t:\tau\mapsto t(\tau)$ is not locally
one-to-one at $\tau=\tau_0$ or $A(\tau_0)=\infty$, then $t_0$ is a
singularity of the differential equation. Thus, singular points
$t_0=t(\tau_0)$ of the differential equation $LA(t)=0$ can occur at
four types of points where
\begin{enumerate}
\item $t_0$ corresponds to an elliptic point of the modular
  curve $X(\Gamma)$, or
\item $t_0$ corresponds to a cusp of $X(\Gamma)$, or
\item $A(t_0)=\infty$, or
\item $t_0$ is a branch point of the covering map
  $X(\Gamma)\to\P^1(\C)$ given by $h(\tau)\to t(\tau)$, where
  $h(\tau)$ is a Hauptmodul of $\Gamma$.
\end{enumerate}
\medskip

\noindent{\bf Examples.}
\begin{enumerate}
\item Ap\'ery's differential equation \eqref{Apery's differential
  equation} has a modular-function parameterization given by
  \begin{equation} \label{Apery's t and A}
    t(\tau)=\left(\frac{\eta(\tau)\eta(6\tau)}{\eta(2\tau)\eta(3\tau)}
    \right)^{12}, \qquad
    A(\tau)=\frac{(\eta(2\tau)\eta(3\tau))^7}
    {(\eta(\tau)\eta(6\tau))^5}.
  \end{equation}
  They are modular on $\Gamma_0(6)+\omega_6$. There are four
  singularities $0,17\pm 12\sqrt2$, and $\infty$. The points $t=0$ and
  $\infty$ correspond to the two cusps $\infty$ and $1/2$. The points
  $t=17-12\sqrt 2$ and $17+12\sqrt 2$ correspond to the elliptic
  points $\tau=i/\sqrt 6$ and $\tau=2/5+i/5\sqrt 6$ of order $2$ fixed
  by $\left(\begin{smallmatrix}0&-1\\6&0\end{smallmatrix}\right)$ and
  $\left(\begin{smallmatrix}12&-5\\30&-12\end{smallmatrix}\right)$,
  respectively.
\item Ap\'ery's differential equation \eqref{Apery's differential
  equation 2} is satisfied by
  \begin{equation} \label{equation: zeta(2) t}
    t(\tau)=q\prod_{n=1}^\infty(1-q^n)^{5\left(\frac n5\right)},
  \end{equation}
  \begin{equation} \label{equation: zeta(2) A}
    A(\tau)=1+\sum_{n=1}^\infty\left(
    \frac{3q^{5n-4}}{1-q^{5n-4}}+\frac{q^{5n-3}}{1-q^{5n-3}}
   -\frac{q^{5n-2}}{1-q^{5n-2}}-\frac{3q^{5n-1}}{1-q^{5n-1}}\right),
  \end{equation}
  where $q=e^{2\pi i\tau}$ and $\left(\frac n5\right)$ is the Legendre
  symbol. The functions are modular on $\Gamma_1(5)$. The
  singularities $0,\infty,(11-5\sqrt 5)/2,(11+5\sqrt 5)/2$ are the
  values of $t(\tau)$ at the four cusps $\infty$, $2/5$, $0$, and
  $1/2$, respectively.
\end{enumerate}

We now briefly review Beukers' argument \cite{Beukers}.
Assume that $t_0=0,~t_1,\dots, t_m$ are singularities of the
differential equation $LA=0$ with $|t_1|<|t_2|<\dots$. Let
$A(t)=1+\sum a_nt^n$ be the $t$-expansion of $A$. In general, the
radius of convergence of the series is $|t_1|$. If the inhomogeneous
differential equation $LB(t)=t$ has a holomorphic solution
$B(t)=t+\sum_{n=2}^\infty b_nt^n$ near $t=0$, this series $B(t)$
in general also has a radius of convergence equal to $|t_1|$. Now
suppose that there is a constant $c$ such that the series $B(t)-cA(t)$
has a larger radius of convergence, i.e., such that $t_1$ is no longer
a singularity of $B(t)-cA(t)$. This would mean that the sequence
$\{b_n/a_n\}$ converges to $c$. Moreover, the larger the ratio
$|t_2/t_1|$ is, the better the rate of convergence is. For example, in
Ap\'ery's differential equation for $\zeta(3)$, we have
$t_1=17-12\sqrt 2$, $t_2=17+12\sqrt2$ and
$t_2/t_1=(17+12\sqrt2)^2=1153.999\dots$. Using the theory of modular
forms, Beukers \cite{Beukers} showed that with the choice of constant
$c=\zeta(3)/6$, the series $B(t)-cA(t)$ no longer has a singularity at
$t_1$. Then the exceptionally large ratio $t_2/t_1$ is sufficient to
imply the irrationality of $\zeta(3)$. We now explain why the constant
$c$ is $\zeta(3)/6$.

\begin{Lemma} \label{lemma: B(t)} Let $\Gamma$ be a discrete subgroup
  of $SL(2,\R)$ commensurable with $SL(2,\Z)$. Let $A(\tau)$ be a
  modular form of weight $k$ and $t(\tau)$ be a non-constant modular
  function on $\Gamma$ such that $t(i\infty)=0$. Let
  $$
    L:\theta^{k+1}+r_k(t)\theta^k+\dots+r_0(t)
  $$
  be the differential operator annihilating $A$. Assume that $g(t)$ is
  a rational function of $t$. Then a solution of the inhomogeneous
  differential equation $LB(t)=g(t)$ is given by
  \begin{equation} \label{equation: B(t)}
    B=A\int^q\left(\cdots\left(\int^q\left(\frac{qdt/dq}t\right)^{k+1}
    \frac{g(t)}A\frac{dq}q\right)\cdots\right)\frac{dq}q,
  \end{equation}
  where the integration is iterated $k+1$ times, and $q=e^{2\pi i\tau}$.
\end{Lemma}

\noindent{\bf Remark.} We did not specify the starting points and the
  paths of integration in \eqref{equation: B(t)} because different
  choices just give different coefficients $d_i$ in general solutions
  $d_0A+d_1\tau A+\cdots+d_k\tau^k A+B_0(t)$, where $B_0(t)$ is a
  fixed solution of $LB(t)=g(t)$. However, when the integrand is a
  holomorphic modular form $\sum_{n=1}^\infty c_nq^n$ that vanishes at
  $i\infty$, we specify the solution to be
  $$
    B=A\int_0^q\left(\cdots
    \left(\int_0^q\left(\sum_{n=1}^\infty c_nq^n\right)
    \frac{dq}q\right)\cdots\right)\frac{dq}q
     =A\sum_{n=1}^\infty\frac{c_n}{n^{k+1}}q^n.
  $$
  In particular, this is what we refer to in the examples in the next
  two sections.

\begin{proof} The proof uses the standard method of variation. Here we
  only prove the case $k=1$; general cases can be proved in the same
  way.

  For convenience, we set
  $$
    G_1=\frac{qdt/dq}t, \qquad G_2=\frac{qdA/dq}A.
  $$
  Then we have
  \begin{equation} \label{equation: lemma B(t)}
    \theta A=A\frac{G_2}{G_1}, \qquad
    \theta\tau=\frac1{2\pi iG_1}.
  \end{equation}
  We will look for two functions $p_1(t)$ and $p_2(t)$ such that
  $B=p_1A+p_2\tau A$ will solve the differential equation
  $LB(t)=g(t)$. We have
  $\theta B=A\theta p_1+\tau A\theta p_2+p_1\theta A+p_2\theta(\tau A)$.
  At this point, we make an additional assumption that $p_1$ and $p_2$
  satisfy
  \begin{equation*} 
    A\theta p_1+\tau A\theta p_2=0
  \end{equation*}
  so that
  $$
    \theta B=p_1\theta A+p_2\theta(\tau A).
  $$
  Differentiating again, we obtain
  $$
    \theta^2B=\theta p_1\theta A+\theta p_2\theta(\tau A)
    +p_1\theta^2A+p_2\theta^2(\tau A).
  $$
  Since $A$ and $\tau A$ are solutions of $L$, we find
  $$
    LB=\theta^2B+r_1\theta B+r_2B=\theta p_1\theta A+\theta
    p_2\theta(\tau A).
  $$
  Therefore, if $p_1$ and $p_2$ satisfy
  $$
    \begin{cases} A\theta p_1+\tau A\theta p_2=0, \\
     \theta A\theta p_1+\theta(\tau A)\theta p_2=g(t), \end{cases}
  $$
  then $B=p_1A+p_2\tau A$ is a solution of $LB(t)=g(t)$. Solving the
  linear equations and using the expressions in \eqref{equation: lemma
  B(t)}, we find
  $$
    \theta p_1=-\frac{2\pi i\tau gG_1}A, \qquad
    \theta p_2=\frac{2\pi igG_1}A
  $$
  It follows that
  $$
    B(t)=-2\pi iA\int^t\frac{\tau gG_1}A\frac{dt}t
    +2\pi i\tau A\int^t\frac{gG_1}A\frac{dt}t.
  $$
  Making the change of variable $t\mapsto q$, we obtain
  $$
    B(t)=-2\pi iA\int^q\frac{\tau gG_1^2}A\frac{dq}q
    +2\pi i\tau A\int^q\frac{gG_1^2}A\frac{dq}q.
  $$
  (Again, although $t\to q$ is a many-to-one mapping, different choices
  of branches just give different solutions.)
  Finally, applying integration by parts to the first integral, we
  conclude that
  $$
    B(t)=A\int^q\left(\int^q\frac{gG_1^2}A\frac{dq}q\right)\frac{dq}q.
  $$
  This proves the lemma for the case $k=1$. General cases can be
  proved in the same way.
\end{proof}

\begin{Lemma} Let all the notations be given as in Lemma \ref{lemma:
  B(t)}. The function
$$
  \left(\frac{qdt/dq}t\right)^{k+1}\frac{g(t)}A
$$
in the above lemma is a modular form of weight $k+2$ on $\Gamma$ with
character $\overline\chi$.
\end{Lemma}

\begin{proof} The lemma is an immediate consequence of the well-known
  property that $(qdt/dq)/t$ is a modular form of weight $2$ on $\Gamma$
  with trivial character.
\end{proof}

\noindent{\bf Example.} Let $t(\tau)$ and $A(\tau)$ be given as in
  \eqref{Apery's t and A}, which satisfy Ap\'ery's differential
  equation \eqref{Apery's differential equation} for $\zeta(3)$. Then
  we have
  $$
    \frac{qdt/dq}t=\frac12\left(
    E_2(\tau)-2E_2(2\tau)-3E_2(3\tau)+6E_2(6\tau)\right),
  $$
  where $E_2(\tau)=1-24\sum_{n=1}^\infty nq^n/(1-q^n)$. By Lemma
  \ref{lemma: B(t)}, a solution of the inhomogeneous solution $LB=t$
  near $t=0$ is given by
  $$
    B(t)=A(t)\int_0^q\int_0^q\int_0^q f(\tau)\frac{dq}q\frac{dq}q
    \frac{dq}q,
  $$
  where
  \begin{equation*}
  \begin{split}
    f(\tau)&=\left(\frac{qdt/dq}t\right)^3\frac{t}{(1-34t+t^2)A} \\
  &=\frac1{240}\left(E_4(\tau)-28E_4(2\tau)+63E_4(3\tau)-36E_4(6\tau)
    \right)
  \end{split}
  \end{equation*}
  and $E_4(\tau)=1+240\sum_{n=1}^\infty n^3q^n/(1-q^n)$ is the
  normalized Eisenstein series of weight $4$ on $SL(2,\Z)$.
\medskip

Since the singularity $t=17-12\sqrt 2$ corresponds to the elliptic
point $\tau=i/\sqrt 6$ fixed by
$\left(\begin{smallmatrix}0&-1\\6&0\end{smallmatrix}\right)$, the
assertion that $B(\tau)-cA(\tau)$, as a function of $t$, is not
singular at $t=17-12\sqrt 2$ is equivalent to that
$B(\tau)-cA(\tau)$ is invariant under the substitution
$\tau\to-1/6\tau$. Now we have $A(-1/6\tau)=-6\tau^2A(\tau)$. Thus,
to prove the latter assertion, it suffices to show that the function
$$
  E(\tau)=\int_0^q\int_0^q\int_0^qf(\tau)\frac{dq}q\frac{dq}q\frac{dq}q
$$
satisfies
\begin{equation} \label{equation: E}
  E(-1/6\tau)-\frac{\zeta(3)}6=-\frac1{6\tau^2}\left(E(\tau)
  -\frac{\zeta(3)}6\right).
\end{equation}
For this purpose, Beukers invoked the standard method of Mellin
transforms.

By the Mellin inversion formula, one can write
$$
  E(-1/6\tau)=\frac1{2\pi i}\int_{2-i\infty}^{2+i\infty}
  \Gamma(s)L(s+3)\left(\frac{2\pi i}{6\tau}\right)^{-s}\,ds,
$$
where $\Gamma(s)$ is the Gamma function and
$$
  L(s)=\zeta(s)\zeta(s-3)(1-28\cdot 2^{-s}+63\cdot 3^{-s}
  -36\cdot 6^{-s})
$$
is the $L$-function associated with the modular form
$f(\tau)$. We then check that the integrand has simple poles only at
$s=0$ and $s=-2$. Moving the line of integration to the left of $\Re
s=-2$, counting residues, making a change of variable $s\mapsto-2-s$,
using the functional equation
$$
  \left(\frac{2\pi}6\right)^{-s}\Gamma(s)L(s)=
 -\left(\frac{2\pi}6\right)^{s-4}\Gamma(4-s)L(4-s),
$$
and then applying the Mellin inversion formula again, we deduce that
\eqref{equation: E} indeed holds. (See \cite[Proposition 1]{Beukers}
for more details. See also the examples in Sections \ref{section:
  elliptic} and \ref{section: cusp}.) This
basically summarizes Beukers' argument.

As we have remarked earlier, the focus of \cite{Beukers} is to give a
modular-function proof of Ap\'ery's irrationality results. Therefore,
the differential equations considered there are very well-behaved in
some sense. That is, the Ap\'ery limits $c$ of them all have the
special property that the the functions $B(t)-cA(t)$ no longer have
singularities at $t_1$. However, in many cases, we do not need such
a strong property in order for the Ap\'ery limits to exist. For
example, consider the differential equation
$$
  \theta^2A+3t(18\theta^2+18\theta+7)A+729t^2(\theta+1)^2A=0.
$$
It has three singularities $0$, $-1/27$, and $\infty$. The local exponents
at these points are $\{0,0\},~\{-2/3,-1/3\}$, and $\{1,1\}$,
respectively. Thus, near 
$t=-1/27$, the solution $A(t)=1+\sum a_nt^n$ has a series expansion
$$
  A(t)=c_1(t+1/27)^{-2/3}+c_2(t+1/27)^{-1/3}+\cdots
$$
for some constants $c_1$ and $c_2$, at least one of which is nonzero.
On the other hand, we can show that the solution $B(t)=t+\sum b_nt^n$ of
$$
  \theta^2B+3t(18\theta^2+18\theta+7)B+729t^2(\theta+1)^2B=t
$$
behaves asymptotically as
$$
  A(t)(d_0+d_1(t+1/27)^{1/3}+\cdots)
$$
near $t=-27$. Therefore, we have
$$
  B(t)-d_0A(t)=A(t)(d_1(t+1/27)^{1/3}+\cdots).
$$
We then can apply the Tauberian theorems to conclude that
$$
  \left|\frac{b_n}{a_n}-d_0\right|\ll\frac1{n^{1/3}}.
$$
To determine the exact value of $d_0$,
we use the theory of modular forms. The detailed
computation is carried out in the next section.

In general, if the singular point $t_1=t(\tau_1)$ closest to the origin
corresponds to an elliptic point $\tau_1$ of order $n$, then the local
exponents at $t_1$ take the form $a_1/n<\dots<a_{k+1}/n$ for some
integers $a_i$, at least one of which is relatively prime to $n$. If
the $(k+1)$-times iterated integral $E(\tau)$ in Lemma \ref{lemma:
B(t)} is a holomorphic function of $\tau$ at $\tau_1$, then the same
argument as that in the previous paragraph will imply that the Ap\'ery
limit is equal to $E(\tau_1)$, provided that the differential equation
does not have another singularity $t_2$ of the same modulus as $t_1$.
Likewise, if the singular point $t_1$ corresponds to a cusp, then the
local exponents are all equal to a rational number $a$. In this case,
we have
$$
  A(t)=(t-t_1)^a\left(c_k\log^k(t-t_1)+c_{k-1}\log^{k-1}(t-t_1)+\cdots
  \right).
$$
Under the same condition on $E(\tau)$, we will also be able to
determine the Ap\'ery limit, although the rate of convergence, in
general, is $|b_n/a_n-d_0|\ll 1/\log n$, which is extremely slow.

In the following two sections, we will determine the Ap\'ery limits of
several examples using the above ideas.

\end{section}

\begin{section}{Elliptic point cases} \label{section: elliptic}

In this section, we will explain in more detail how to determine the
Ap\'ery limits when the singularity $t_1$ closest to the origin
corresponds to an elliptic point $\tau_1$. Throughout the section, we
assume that $\Gamma$ is a discrete subgroup of $SL(2,\R)$ of genus
zero commensurable with $SL(2,\Z)$, $A(\tau)$ is a (meromorphic)
modular form of weight $k$ with character $\chi$ on $\Gamma$, and
$t(\tau)$ is a uniformizer of the modular curve $X(\Gamma)$ so that
the differential equation $LA(t)=0$ satisfied by $A$ and $t$ has
polynomial coefficients. We also assume that $A(i\infty)=1$ and
$t(i\infty)=0$ so that the differential equation has local exponents
all equal to $0$ at $t=0$. We let $B(t)=t^j+\cdots$ be the solution of
$LB(t)=j^{k+1}t^j$ holomorphic at $t=0$, where $j$ is a positive
integer less than the maximal degree of the coefficients of $L$. For a
function $f$, with a slight abuse of notations, we write $f(\tau)$ if
we consider $f$ as a function of $\tau$, and write $f(t)$ if we
consider it as a (multi-valued) function of $t$.

Let $h(t)$ be the coefficient of $\theta^{k+1}$ in the differential
equation $LA=0$. We assume that the integrand
$$
  f(\tau)=\left(\frac{qdt/dq}t\right)^{k+1}\frac{t^j}{Ah(t)}
$$
inside \eqref{equation: B(t)} is a holomorphic modular form so that
its $L$-function converges in some half-plane. Note that, by the
assumption that $t(i\infty)=0$ and $A(i\infty)=1$, the constant term of
the Fourier expansion $\sum_{n=1}^\infty c_nq^n$ of $f(\tau)$ is $0$.
According to the discussion in the previous section, the Ap\'ery
limit is equal to the value of the $(k+1)$-times iterated integral
$$
  E(\tau)=\int_0^q\left(\cdots\left(\int_0^q f(\tau)
    \frac{dq}q\right)\cdots\right)\frac{dq}q
  =\sum_{n=1}^\infty\frac{c_n}{n^{k+1}}q^n
$$
at $\tau_1$. Now to evaluate $E(\tau_1)$, we express $E(\tau)$ using
the Mellin inversion formula
$$
  e^{-x}=\frac1{2\pi i}\int_{c-i\infty}^{c+i\infty}\Gamma(s)x^{-s}\,ds
$$
which holds for all $c>0$ and all complex numbers $x$ with $\Re x>0$,
and then try to deduce information about $E(\tau_1)$ by complex
analytic argument. For example, if $f(\tau)$ is an eigenfunction with
eigenvalue $\epsilon$ of the Atkin-Lehner involution
$\omega_N=\left(\begin{smallmatrix}0&-1\\N&0\end{smallmatrix}\right)$
for some positive integer $N$, then by an argument similar to that in
the proof of \eqref{equation: E}, we can show that
$$
  E(-1/N\tau)=\epsilon^{-1}(\sqrt N\tau)^{-k}E(\tau)+(\text{residues}).
$$
Then we set $\tau=i/\sqrt N$ to get the value of $E(i/\sqrt N)$
(provided that $i^k\epsilon\neq 1$).

When the elliptic element that fixes $\tau_1$ is not
$\left(\begin{smallmatrix}0&-1\\N&0\end{smallmatrix}\right)$, the
situation is more complicated. Assume that $f(\tau)$ is modular on
$\Gamma_0(N)$. If we apply the Mellin inversion formula in a
straightforward manner and write
$$
  E(\tau)=\frac1{2\pi i}\int_{c-i\infty}^{c+i\infty}
  \Gamma(s)L(s+k+1,f)\left(\frac{2\pi\tau}i\right)^{-s}\,ds,
$$
the complex analytic argument will yield
$$
  E(\tau)=\frac\epsilon{\tau^k}\sum_{n=1}^\infty\frac{c^\ast_n}{n^{k+1}}
  e^{-2\pi in/N\tau}+(\text{residues}),
$$
where $c_n^\ast$ denote the Fourier coefficients of $(\sqrt
N\tau)^{-k-2}f(-1/N\tau)$ (which is a modular form on $\Gamma_0(N)$
since $\omega_N$ normalizes $\Gamma_0(N)$). At this point, it is not
clear how one should proceed to obtain information about $E(\tau_1)$.
In such a situation, we need to apply the Mellin inversion formula in a
different way. The key observation is the following simple fact.

\begin{Lemma} \label{lemma: alternative} For a matrix
  $\left(\begin{smallmatrix}a&b\\c&d\end{smallmatrix}\right)\in
  SL(2,\R)$ with $c>0$ and for a pair of functions $g,g^\ast:\H\to\C$,
  the relation
$$
  g\left(\frac{a\tau+b}{c\tau+d}\right)=\epsilon(c\tau+d)^k g^\ast(\tau)
$$
holds for some constant $\epsilon$ and some integer $k$ if and only if
$$
  g\left(\frac{\tau}c+\frac ac\right)
 =\epsilon(-\tau)^{-k}g^\ast\left(-\frac1{c\tau}-\frac dc\right)
$$
holds.
\end{Lemma}

\begin{proof} We have
$$
  \frac\tau c+\frac ac=\begin{pmatrix}a&b\\c&d\end{pmatrix}
  \left(-\frac1{c\tau}-\frac dc\right).
$$
This proves the lemma.
\end{proof}

Now if $g(\tau)=\sum_{n=0}^\infty c_nq^n$ and
$g^\ast(\tau)=\sum_{n=0}^\infty c_n^\ast q^n$ in the above lemma are
holomorphic modular forms, then the lemma yields a functional equation
between the pair of functions
$$
  L(s)=\sum_{n=1}^\infty\frac{c_n}{n^s}e^{2\pi ina/c}, \qquad
  L^\ast(s)=\sum_{n=1}^\infty\frac{c_n^\ast}{n^s}e^{-2\pi ind/c}
$$
(See Lemma \ref{lemma: functional equation} below.)
Using the functional equation, we deduce that
$$
  E\left(-\frac1{c\tau}-\frac dc\right)=\frac\epsilon{\tau^k}
  E\left(\frac\tau c+\frac ac\right)+(\text{residues}).
$$
We then make a suitable choice of $\tau$ such that
$-1/c\tau-d/c=\tau/c+a/c=\tau_1$ and an evaluation of $E(\tau_1)$
follows.

We now prove the functional equation between $L(s)$ and $L^\ast(s)$ in
the following lemma. The lemma may have appeared somewhere in literature.
However, failing to locate such a reference, we give a complete
proof here. Notice that when $g=g^\ast$ is a modular form of weight $k$ on
$\Gamma_0(N)$ that is an eigenfunction of the Atkin-Lehner involution
$w_N$ with eigenvalue $\epsilon$, the lemma gives the familiar
functional equation
$$
  i^k\epsilon\left(\frac{2\pi}{\sqrt N}\right)^{-s}\Gamma(s)L(s,g)
 =\left(\frac{2\pi}{\sqrt N}\right)^{s-k}\Gamma(k-s)L(k-s,g)
$$
for the $L$-function $L(s,g)$ of $g$.

\begin{Lemma} \label{lemma: functional equation} Let $\Gamma$ and
  $\Gamma^\ast$ be two discrete subgroups of $SL(2,\R)$ commensurable
  with $SL(2,\Z)$ such that the cusp $\infty$ has width $1$. Assume
  that $g(\tau)=\sum_{n=0}^\infty c_nq^n$ and
  $g^\ast(\tau)=\sum_{n=0}^\infty c_n^\ast q^n$ are holomorphic
  modular forms of weight $k$ on $\Gamma$ and $\Gamma^\ast$,
  respectively. Suppose that
  $\gamma=\left(\begin{smallmatrix}a&b\\c&d\end{smallmatrix}
  \right)$ with $c>0$ is a matrix in $SL(2,\R)$ such that
  \begin{equation} \label{equation: lemma functional equation}
    g\left(\frac{a\tau+b}{c\tau+d}\right)=\epsilon(c\tau+d)^k
    g^\ast(\tau)
  \end{equation}
  for some constant $\epsilon$. Then the two functions
  $$
    L(s)=\sum_{n=1}^\infty\frac{c_n}{n^s}e^{2\pi ina/c}, \qquad
    L^\ast(s)=\sum_{n=1}^\infty\frac{c_n^\ast}{n^s}e^{-2\pi ind/c}.
  $$
  satisfy the functional equation
  \begin{equation} \label{equation: functional equation}
    \left(\frac{2\pi}c\right)^{-s}\Gamma(s)L(s)
   =i^k\epsilon\left(\frac{2\pi}c\right)^{s-k}\Gamma(k-s)L^\ast(k-s).
  \end{equation}
\end{Lemma}

\begin{proof} For the ease of representation, here we only prove the
  case where the constant terms $c_0$ and $c_0^\ast$ are $0$. The case
  $c_0,c_0^\ast\neq 0$ can be proved by a minor modification.

First of all, setting $\tau=(-d\tau+b)/(c\tau-a)$ in \eqref{equation:
  lemma functional equation}, we obtain
$$
  g^\ast\left(\frac{-d\tau+b}{c\tau-a}\right)
  =\epsilon^{-1}(-1)^k(c\tau-a)^k g(\tau),
$$
which by Lemma \ref{lemma: alternative} is equivalent to
\begin{equation} \label{equation: lemma 4 temp}
  g^\ast\left(\frac\tau c-\frac dc\right)
  =\epsilon^{-1}\tau^{-k} g\left(-\frac1{c\tau}+\frac ac\right).
\end{equation}
We now consider the integral
$$
  \int_0^\infty y^{s-1}g\left(\frac{iy}c+\frac ac\right)\,dy.
$$
We have
\begin{equation*}
\begin{split}
  \int_0^\infty y^{s-1}g\left(\frac{iy}c+\frac ac\right)\,dy
  &=\sum_{n=1}^\infty c_ne^{2\pi ina/c}\int_0^\infty y^{s-1}
    e^{-2\pi ny/c}\,dy \\
  &=\left(\frac{2\pi}c\right)^{-s}\Gamma(s)
    \sum_{n=1}^\infty\frac{c_n}{n^s}e^{2\pi ina/c}
   =\left(\frac{2\pi}c\right)^{-s}\Gamma(s)L(s).
\end{split}
\end{equation*}
Now we break the integral into two parts as usual, one from $0$ to
$1$ and the other from $1$ to $\infty$. For the integral from $0$
to $1$, we make a change of variable $y\mapsto 1/y$ and obtain
$$
  \int_0^{1}y^{s-1}g\left(\frac{iy}c+\frac ac\right)\,dy
 =\int_{1}^\infty y^{1-s}g\left(\frac{i}{cy}+\frac ac\right)\frac{dy}{y^2}.
$$
Substituting the relation \eqref{equation: lemma 4 temp} with
$\tau=iy$ into the last expression, we see that
$$
  \int_0^{1}y^{s-1}g\left(\frac{iy}c+\frac ac\right)\,dy
 =i^k\epsilon\int^\infty_{1}y^{k-s-1}g^\ast\left(\frac{iy}c
  -\frac dc\right)\,dy
$$
and
$$
  \left(\frac{2\pi}c\right)^{-s}\Gamma(s)L(s)
 =\int^\infty_{1}y^{s-1}g\left(\frac{iy}c+\frac ac\right)\,dy
  +i^k\epsilon \int^\infty_{1}
  y^{k-s-1}g^\ast\left(\frac{iy}c-\frac dc\right)\,dy.
$$
By the same token, we can also show that
\begin{equation*}
\begin{split}
  \left(\frac{2\pi}c\right)^{-s}\Gamma(s)L^\ast(s)
&=\int^\infty_{1}y^{s-1}g^\ast\left(\frac{iy}c-\frac dc\right)\,dy \\
&\qquad+\epsilon^{-1}(-i)^k\int^\infty_{1}
  y^{k-s-1}g\left(\frac{iy}c+\frac ac\right)\,dy.
\end{split}
\end{equation*}
Upon the substitution $s\mapsto k-s$ into the last expression, we
immediately get the claimed functional equation. This completes
the proof.
\end{proof}

We now give two examples to illustrate our method.
\medskip

\noindent{\bf Example 1.} Consider the differential equation
\begin{equation} \label{equation: case h}
  \theta^2A+3t(18\theta^2+18\theta+7)A+729t^2(\theta+1)^2A=0.
\end{equation}
This is case ($h$) that \cite{ASZ} fails to identify the Ap\'ery
limit. According to \cite{ASZ}, Arne Meurman conjectures that the limit
is $2\pi^2/81-L(2,\chi_3)/2$, where $L(s,\chi_3)$ is the Dirichlet
$L$-function associates with $\chi_3=\left(\frac\cdot 3\right)$.
Meurman also suggests the following Ramanujan-like formula
$$
  \sum_{n=1}^\infty\left(\frac n3\right)\frac
  {(-e^{-\pi/\sqrt 3})^n}{n^2(1-(-e^{-\pi/\sqrt 3})^n)}
  =\frac{2\pi^2}{81}-\frac12L(2,\chi_3).
$$
Here we prove that these conjectures are indeed correct.

We first note that \eqref{equation: case h} is the differential
equation satisfied by
$$
  t(\tau)=\frac{\eta(3\tau)^{12}}{\eta(\tau)^{12}}, \qquad
  A(\tau)=\frac1{1+27t(\tau)}
    \left(1+6\sum_{n=1}^\infty\left(\frac n3\right)
    \frac{q^n}{1-q^n}\right),
$$
where $A(\tau)$ is a modular form of weight $1$ on $\Gamma_0(3)$ with
character
\begin{equation} \label{equation: character}
  \chi\begin{pmatrix}a&b\\c&d\end{pmatrix}
 =\left(\frac d3\right).
\end{equation}
The singular points $0$, $\infty$, and $-1/27$ correspond to the cusps
$\infty$, $0$, and the elliptic point
$$
  \tau_0=\frac{3+\sqrt{-3}}6,
$$
respectively. Using Lemma \ref{lemma: B(t)} we find a solution of
$$
  \theta^2B+3t(18\theta^2+18\theta+7)B+729t^2(\theta+1)^2B=t
$$
is given by
$$
  B(t)=A(t)\int_0^q\int_0^q\left(\frac{qdt/dq}t\right)^2\frac t
  {A(t)(1+27t)^2}\frac{dq}q\frac{dq}q
  =A(t)\int_0^q\int_0^q\frac{\eta(3\tau)^9}{\eta(\tau)^3}\frac{dq}q
  \frac{dq}q.
$$
The local exponents of the differential equation at $-1/27$ are $-1/3$
and $-2/3$. Thus, near $t=-1/27$, we have the series expansion
$$
  A(t)=c_1(t+1/27)^{-2/3}+c_2(t+1/27)^{-1/3}+\cdots
$$
for some constants $c_i$ with at least one of $c_1$ and $c_2$ being
nonzero. Now write
$$
  E(\tau)=\int_0^q\int_0^q\frac{\eta(3\tau)^9}{\eta(\tau)^3}\frac{dq}q
  \frac{dq}q.
$$
Regardless of what the constants $c_1$ and $c_2$ are, we have
$$
  B(t)-E(\tau_0)A(t)=A(t)(E(t)-E(\tau_0))=A(t)(0+d_1(t+1/27)^{1/3}
  +\cdots).
$$
At this point, the Tauberian theorems (using, for example,
\cite[Corollary 1.7.3]{BGT}) already assert that if $A(t)=1+\sum_n
a_nt^n$ and $B(t)=t+\sum_n b_nt^n$, then we have
$$
  \left|\frac{\sum_{n=1}^N b_n}{\sum_{n=1}^N a_n}-E(\tau_0)\right|
  \ll N^{-1/3}
$$
as $N\to\infty$. To get the stronger statement $b_n/a_n\to E(\tau_0)$,
we consider the sequences $a_n'=a_n-a_{n-1}$ and $b_n'=b_n-b_{n-1}$,
whose generating functions are $(1-t)A(t)$ and $(1-t)B(t)$,
respectively. Since these two functions have similar asymptotic
behaviors as $A(t)$ and $B(t)$ near $t=-1/27$, the same Tauberian
theorem implies that
$$
  \left|\frac{b_n}{a_n}-E(\tau_0)\right|\ll n^{-1/3}.
$$
In particular, the Ap\'ery limit is $E(\tau_0)$. It remains to
determine the value of $E(\tau_0)=E((3+\sqrt{-3})/6)$.

The integrand $f(\tau)=\eta(3\tau)^9/\eta(\tau)^3$ in the definition
of $E(\tau)$ is in fact the Eisenstein series of weight $3$ associated
with the cusp $0$ with the same character as \eqref{equation:
character}. It can be written as
$$
  \sum_{n=1}^\infty\frac{n^2(q^n-q^{2n})}{1-q^{3n}}
 =\sum_{n=1}^\infty n^2\sum_{k=1}^\infty\left(\frac k3\right)q^{kn}.
$$
It follows that
$$
  E(\tau)=\sum_{n=1}^\infty\sum_{k=1}^\infty\frac1{k^2}
  \left(\frac k3\right)q^{kn}.
$$
Our strategy is to use Lemma \ref{lemma: functional equation} to show
that
$$
  E(-1/3\tau+1/3)=\frac1\tau E(\tau/3+2/3)+r(\tau)
$$
for some function $r(\tau)$ involving $\tau$ and special values of
Dirichlet $L$-functions. Then we set $\tau=(-1+\sqrt{-3})/2$. With
this choice of $\tau$, we have
$$
  -\frac1{3\tau}+\frac13=\frac\tau3+\frac23=\frac{3+\sqrt{-3}}6=\tau_0.
$$
From this we can obtain the values of $E(\tau_0)$.

Write the Fourier coefficients of $q^n$ in $f(\tau)$ as $c_n$.
We apply Lemma \ref{lemma: functional equation} with $g=g^\ast=f$,
$$
  \gamma=\begin{pmatrix}a&b\\c&d\end{pmatrix}
 =\begin{pmatrix}2&-1\\3&-1\end{pmatrix},
$$
$\epsilon=g(\chi)=-1$, and define
$$
  L(s)=\sum_{n=1}^\infty\frac{c_ne^{4\pi in/3}}{n^s}, \qquad
  L^\ast(s)=\sum_{n=1}^\infty\frac{c_ne^{2\pi in/3}}{n^s}.
$$
By the Mellin inversion formula, we have
\begin{equation} \label{equation: Mellin}
\begin{split}
  E(-1/3\tau+1/3)&=\sum_{n=1}^\infty\frac{c_ne^{2\pi in/3}}{n^2}
  e^{-2\pi in/3\tau} \\
  &=\frac1{2\pi i}\int_{3/2-i\infty}^{3/2+i\infty}
    \Gamma(s)L^\ast(s+2)\left(\frac{2\pi i}{3\tau}\right)^{-s}\,ds.
\end{split}
\end{equation}
We then move the line of integration to $\Re s=-5/2$ and make a change
of variable $s\mapsto -1-s$. (By the expression of $L(s)$ in terms of
$\zeta(s)$ and $L(s,\chi_3)$ given later in \eqref{equation: L(s)}, we
see that this is justified.) We obtain
$$
  E(-1/3\tau+1/3)=(\text{residues})+
  \frac1{2\pi i}\int_{3/2-i\infty}^{3/2+i\infty}
  \Gamma(-1-s)L^\ast(1-s)\left(\frac{2\pi i}{3\tau}\right)^{s+1}\,ds.
$$
Now Lemma \ref{lemma: functional equation} implies that
$$
  \left(\frac{2\pi}3\right)^{-s}\Gamma(s)L^\ast(s)
 =-i\left(\frac{2\pi}3\right)^{s-3}\Gamma(3-s)L(3-s).
$$
It follows that
\begin{equation*}
\begin{split}
  \Gamma(-1-s)L^\ast(1-s)&=\frac{\Gamma(1-s)L^\ast(1-s)}{s(s+1)}
   =-i\left(\frac{2\pi}3\right)^{-1-2s}\frac{\Gamma(s+2)L(s+2)}
   {s(s+1)} \\
  &=-i\left(\frac{2\pi}3\right)^{-1-2s}\Gamma(s)L(s+2),
\end{split}
\end{equation*}
and we have
$$
  E(-1/3\tau+1/3)=(\text{residues})+
  \frac1{2\pi i\tau}\int_{3/2-i\infty}^{3/2+i\infty}
  \Gamma(s)L(s+2)\left(\frac{2\pi\tau}{3i}\right)^{-s} \,ds.
$$
By the Mellin inversion formula again, the integral in the last
expression is precisely $E(\tau/3+2/3)$. That is, we have
\begin{equation} \label{equation: lemma temp1}
  E(-1/3\tau+1/3)=(\text{residues})+\frac1\tau E(\tau/3+2/3).
\end{equation}
It remains to compute the residues.

We have
\begin{equation*}
\begin{split}
  L^\ast(s)=\sum_{n=1}^\infty\sum_{k=1}^\infty n^2\left(\frac k3\right)
  \frac{e^{2\pi ink/3}}{n^sk^s}.
\end{split}
\end{equation*}
Partitioning the double sum according to the residue classes modulo
$3$ and simplifying, we obtain
\begin{equation} \label{equation: L(s)}
  L^\ast(s)=-\frac12\zeta(s-2)(1-3^{3-s})L(s,\chi_3)+\frac{\sqrt3i}2
  L(s-2,\chi_3)\zeta(s)(1-3^{-s}).
\end{equation}
Using the fact that $\zeta(s)$ has zeros at negative even integers and
$L(s,\chi_3)$ has zeros at negative odd integers, we see that the
integrand
\begin{equation*}
\begin{split}
 &\Gamma(s)L^\ast(s+2)\left(\frac{2\pi i}{3\tau}\right)^{-s}
 =\Gamma(s)\left(\frac{2\pi i}{3\tau}\right)^{-s} \\
 &\qquad\qquad\times
  \left(-\frac12\zeta(s)(1-3^{1-s})L(s+2,\chi_3)+\frac{\sqrt3i}2
  L(s,\chi_3)\zeta(s+2)(1-3^{-2-s})\right)
\end{split}
\end{equation*}
in \eqref{equation: Mellin} has poles only at $s=0$ and $s=-1$. We
find that the residue at $s=0$ is
\begin{equation} \label{equation: residue 1}
  \operatorname{Res}_{s=0}
 =-\frac12\zeta(0)(1-3)L(2,\chi_3)+\frac{\sqrt 3i}2L(0,\chi_3)
   \zeta(2)(1-1/9).
\end{equation}
With the evaluations $\zeta(0)=-1/2$, $L(0,\chi_3)=1/3$, and
$\zeta(2)=\pi^2/6$, we simplify the expression into
$$
  -\frac12L(2,\chi_3)+\frac{2\sqrt 3\pi^2i}{81}.
$$
Also, the residue at $s=-1$ is
$$
  -\left(-\frac12\zeta(-1)(1-9)L(1,\chi_3)+\frac{\sqrt 3i}2L'(-1,\chi_3)
  (1-1/3)\right)\frac{2\pi i}{3\tau}.
$$
Using the functional equation
$$
  \left(\frac\pi3\right)^{-(s+1)/2}\Gamma\left(\frac{s+1}2\right)
  L(s,\chi_3)=\left(\frac\pi 3\right)^{-(2-s)/2}\Gamma
  \left(\frac{2-s}2\right)L(1-s,\chi_3)
$$
for $L(s,\chi_3)$, one may deduce that
$$
  L'(-1,\chi_3)=\frac{3\sqrt 3}{4\pi}L(2,\chi).
$$
Together with the evaluations $\zeta(-1)=-1/12$ and
$L(1,\chi_3)=\pi\sqrt 3/9$, we find that the residue at $s=-1$ can be
simplified as
\begin{equation} \label{equation: residue 2}
  \operatorname{Res}_{s=-1}=
  \frac{2\pi i\sqrt3}{81\tau}+\frac1{2\tau}L(2,\chi).
\end{equation}
In summary, the above computation \eqref{equation: lemma temp1},
\eqref{equation: residue 1}, \eqref{equation: residue 2} shows that
$$
  E(-1/3\tau+1/3)=-\frac{\tau-1}{2\tau}L(2,\chi)
  +\frac{\tau+1}\tau\frac{2\pi^2i\sqrt 3}{81}+\frac1\tau E(\tau/3+2/3).
$$
Setting $\tau=(-1+\sqrt{-3})/2$, we find
$$
  E\left(\frac{3+\sqrt{-3}}6\right)=\frac{2\pi^2}{81}-\frac12L(2,\chi_3).
$$
We summarize the above computation as follows.
\medskip

\begin{Theorem} Let $\{a_n\}$ and $\{b_n\}$ be the
  sequences satisfying the recursive relation
$$
  (n+2)^2u_{n+2}+(54n^2+162n+129)u_{n+1}+(n+1)^2u_n=0 \quad
  (u_n=a_n\text{ or }b_n)
$$
with the initial values $a_0=1$, $a_1=-21$, $b_0=0$, and $b_1=1$. Then
we have
$$
  \left|\frac{b_n}{a_n}-\frac{2\pi^2}{81}+\frac12L(2,\chi_3)\right|
  \ll n^{-1/3}
$$
as $n\to\infty$. Moreover, the evaluation
$$
  \sum_{n=1}^\infty\left(\frac n3\right)\frac
  {(-e^{-\pi/\sqrt 3})^n}{n^2(1-(-e^{-\pi/\sqrt 3})^n)}
  =\frac{2\pi^2}{81}-\frac12L(2,\chi_3)
$$
holds.
\end{Theorem}

\noindent{\bf Example 2.} In this example, we will construct a
differential equation whose Ap\'ery limits give a good rational
approximation to the value $L(2,f)$ at $2$ of the $L$-function
associated with the cusp form $f(\tau)=\eta(\tau)^3\eta(7\tau)^3$.

The Hauptmodul of the genus-zero group $\Gamma_0(7)+\omega_7$ is given
by $\eta(\tau)^4/\eta(7\tau)^4+49\eta(7\tau)^4/\eta(\tau)^4$. We set
$$
  t(\tau)=\left(\frac{\eta(\tau)^4}{\eta(7\tau)^4}+14
  +49\frac{\eta(7\tau)^4}{\eta(\tau)^4}\right)^{-1}
  =q-10q^2+49q^3-184q^4+\cdots.
$$
The value of $t(\tau)$ at the unique cusp $\infty$ is obviously $0$,
and at the elliptic points $i/\sqrt 7$, $(5+\sqrt{-3})/14$, and
$(7+\sqrt{-7})/14$, it takes the values $1/28$, $1$, and $\infty$,
respectively.

According to Lemma \ref{lemma: B(t)} and the discussion at the
beginning of this section, we can basically choose any rational
function $g(t)$ such that $g(0)=0$ and its numerator has a degree
smaller than its denominator and then set
$$
  A(\tau)=\left(\frac{qdt/dq}t\right)^2\frac{g(t)}{f(\tau)}.
$$
Then the differential equation satisfied by $t$ and $A$ will give us a
rational approximation to $L(2,f)$.

Here we choose $g(t)=t/(1-29t+28t^2)$. With this choice, we find
$$
  A(\tau)=1+2\sum_{n=1}^\infty\left(\frac n7\right)\frac{q^n}{1-q^n}
$$
and the differential equation satisfied by $t$ and $A$ is
$$
  (1-t)^2(1-28t)\theta^2A-14t(1-t)^2\theta A-(2t+4t^2)A=0.
$$
The solution $B(t)$ of the inhomogeneous differential equation
$$
  (1-t)^2(1-28t)\theta^2B-14t(1-t)^2\theta B-(2t+4t^2)B=t(1-t)
$$
has a $t$-expansion $B(t)=t+45t^2/4+\cdots$. Furthermore, we can show
that $B(t)-L(2,f)A(t)$ has no singularity at $t=1/28$. Therefore, we
have the following good rational approximation to $L(2,f)$.

\begin{Theorem} Let $f(\tau)=\eta(\tau)^3\eta(7\tau)^3$ be the cusp
  form of weight $3$ with character
  $\chi\left(\begin{smallmatrix}a&b\\c&d\end{smallmatrix}\right)=
   \left(\frac d7\right)$ on $\Gamma_0(7)$ and $L(s,f)$ be its
  associated $L$-function. Let $\{a_n\}$ and
  $\{b_n\}$ be sequences satisfying
$$
  (n+3)^2u_{n+3}-(30n^2+134n+150)u_{n+2}+(57n^2+142n+81)u_{n+1}
  -(14n+28n^2)u_n=0
$$
with the initial values
$$
  a_0=1, \quad a_1=2, \quad a_2=24, \quad b_0=0, \quad b_1=1, \quad
  b_2=45/4.
$$
Then we have
$$
  \left|\frac{b_n}{a_n}-L(2,f)\right|\ll 28^{-n}.
$$
\end{Theorem}
\end{section}

\begin{section}{Cusp cases} \label{section: cusp}
Let $t(\tau)$ and $A(\tau)$ be given and assume that the integrand in
\eqref{equation: B(t)} is a holomorphic modular form that vanishes at
$i\infty$ as before. In this section we consider the cases where the
singularity $t_1$ closest to the origin corresponds to a cusp $\alpha$
of the modular curve $X(\Gamma)$. Pick an element in
$\sigma=\left(\begin{smallmatrix}a&b\\c&d\end{smallmatrix}\right)\in
SL(2,\R)$ (preferably lying in the normalizer of $\Gamma$ in
$SL(2,\R)$) such that $\sigma\infty=\alpha$. Note that the asymptotic
behavior of $A(\tau)$ near $\alpha$ is determined by that of
$A(\sigma\tau)$ in a neighborhood of $\infty$. To be more precise, we
have $A(\sigma\tau)=(c\tau+d)^k A^\ast(\tau)$ for some modular form
$A^\ast(\tau)$ of weight $k$ on the group
$\sigma^{-1}\Gamma\sigma$. Also, $t^\ast(\tau)=t(\sigma\tau)-t_1$ is
an algebraic function of $t$ that vanishes at $\infty$. Then for some
positive integer $r$ depending on the order of $A^\ast(\tau)$ at
$\infty$, $A^\ast(\tau)^r$ has a power series expansion in
$t^\ast$. The term 
$(c\tau+d)^k$ gives a logarithmic factor $d_k\log^k
t^\ast+\cdots+d_0$. Therefore
$$
  A(\sigma\tau)=(d_k\log^k t^\ast+\cdots+d_0)(t^\ast)^a\times(\text{a power
  series in }t^\ast)
$$
for some rational number $a$. (The number $a$ in fact is equal to the
local exponent of the differential equation at $t_1$.) This gives us
the asymptotic behavior of $A(\tau)$ near $\alpha$.

By the same token, to determine the asymptotic behavior of $B(t)$
near $t=t_1$, we need to consider $B(\sigma\tau)$ in a neighborhood of
$\infty$. Let $f(\tau)$ be the modular form of weight $k+2$ inside the
integral in \eqref{equation: B(t)}. We apply Lemma \ref{lemma:
  functional equation} to $f(\tau)$ and
$f^\ast(\tau)=(c\tau+d)^{-k-2}f(\tau)$, which is a modular form on
$\sigma^{-1}\Gamma\sigma$. We then use the functional equation to
find the asymptotics of $B(t)$. The actual computation is similar to
that for the elliptic point cases. We now give some examples.
\medskip

\noindent{\bf Example 3.} Consider the differential equation
$$
  \theta^2A-4t(8\theta^2+8\theta+3)A+256t^2(\theta+1)^2A=0,
$$
which is case ($e$) in \cite{ASZ}. It is conjectured that the Ap\'ery
limit is $L(2,\chi_{-1})/2$, where $\chi_{-1}$ is the odd Dirichlet
character modulo $4$. We now show that this is indeed the case.

The differential equation is the one satisfied by
$$
  t(\tau)=\frac{\eta(\tau)^8\eta(4\tau)^{16}}{\eta(2\tau)^{24}},
  \qquad
  A(\tau)=\frac{\eta(2\tau)^{22}}{\eta(\tau)^{12}\eta(4\tau)^8},
$$
where $A(\tau)$ is a meromorphic modular form of weight $1$ on
$\Gamma_0(4)$ with character
\begin{equation} \label{equation: cusp example character}
  \chi\begin{pmatrix}a&b\\4c&d\end{pmatrix}=(-1)^c\chi_{-1}(d),
\end{equation}
and $t(\tau)$ is a modular function on $\Gamma_0(4)$. By Lemma
\ref{lemma: B(t)},
$$
  B(\tau)=A(\tau)\int_0^q\int_0^q\frac{\eta(\tau)^4\eta(4\tau)^8}
  {\eta(2\tau)^6}\frac{dq}q\frac{dq}q
$$
is a solution of the inhomogeneous
$\theta^2B-4t(8\theta^2+8\theta+3)B+256t^2(\theta+1)^2B=t$.

The singularities $0,~1/16,~\infty$ correspond to the cusps
$\infty,~0$, and $1/2$, respectively, and the local exponents are
$\{0,0\}$, $\{-1/2,-1/2\}$, and $\{1,1\}$. Since
$$
  A(-1/4\tau)=\frac{\tau}{2i}\frac{\eta(2\tau)^{22}}{\eta(\tau)^8
  \eta(4\tau)^{12}},
$$
we know that
$$
  A(t)=(t-1/16)^{-1/2}(c_1\log(t-1/16)+c_2+\cdots)
$$
as $t$ approaches $1/16$ for some nonzero constant $c_1$. (The exact
value of $c_1$ is not important, as it suffices to know that it is not
zero.) Let
$$
  E(\tau)=\int_0^q\int_0^q\frac{\eta(\tau)^4\eta(4\tau)^8}
  {\eta(2\tau)^6}\frac{dq}q\frac{dq}q.
$$
In the following we will show that
\begin{equation} \label{equation: cusp example 1}
  E(-1/4\tau)=\frac12L(2,\chi_{-1})-\frac{\pi^2i}{32\tau}
  +\frac{i}{8\tau}\int_0^q\int_0^q\frac{\eta(\tau)^8\eta(4\tau)^4}
  {\eta(2\tau)^6}\frac{dq}q\frac{dq}q.
\end{equation}
This would imply that
$$
  B(t)-\frac12L(2,\chi_{-1})A(t)=(t-1/16)^{-1/2}(d_0+\cdots).
$$
Then we deduce from the Tauberian theorems that
$$
  \left|\frac{b_n}{a_n}-\frac12L(2,\chi_{-1})\right|\ll\frac1{\log n}.
$$
(Again, we shall apply the Tauberian theorems to the sequence
$a_n'=a_n-a_{n-1}$ and $b_n'=b_n-b_{n-1}$ with generating functions
$(1-t)A(t)$ and $(1-t)B(t)$ to get the stronger conclusion.) We now
prove \eqref{equation: cusp example 1}.

The function $\eta(\tau)^4\eta(4\tau)^8/\eta(2\tau)^6$ is in fact the
Eisenstein series of weight $3$ associated with the cusp $1/2$ with
character given by \eqref{equation: cusp example character}. Its
$q$-expansion can be alternatively written as
$$
  \frac{\eta(\tau)^4\eta(4\tau)^8}{\eta(2\tau)^6}
 =\sum_{n=1}^\infty(-1)^{n-1}\frac{n^2q^n}{1+q^{2n}}.
$$
Then by the Mellin inversion formula
\begin{equation*}
\begin{split}
  E(-1/4\tau)=\frac1{2\pi i}\int_{3/2-i\infty}^{3/2+i\infty}
  \Gamma(s)L(s+2)\left(\frac{2\pi i}{4\tau}\right)^{-s}\,ds,
\end{split}
\end{equation*}
where
$$
  L(s+2)=\sum_{n=1}^\infty(-1)^{n-1}n^2\sum_{k=1}^\infty
  \frac{(-1)^{k-1}}{((2k-1)n)^{s+2}}=\zeta(s)L(s+2,\chi_{-1})(1-2^{1-s}).
$$
Moving the path of integration to $\Re s=-5/2$ and making a change of
variable $s\mapsto -1-s$, we obtain
\begin{equation*}
\begin{split}
 &E(-1/4\tau)=\frac12L(2,\chi_{-1})-\frac{\pi^2i}{32\tau} \\
 &\qquad+\frac1{2\pi i}\int_{3/2-i\infty}^{3/2+i\infty}
  \Gamma(-1-s)\zeta(-1-s)L(1-s,\chi_{-1})(1-2^{2+s})
  \left(\frac{2\pi i}{4\tau}\right)^{s+1}\,ds.
\end{split}
\end{equation*}
Using the functional equations for $\zeta(s)$ and $L(s,\chi_{-1})$ and
the Legendre duplication formula for $\Gamma(s)$, we find the integral
on the right-hand side is equal to
$$
  \frac1{2\pi i}\frac i{2\tau}\int_{3/2-i\infty}^{3/2+i\infty}\Gamma(s)
  \zeta(s+2)L(s,\chi_{-1})(1-2^{-s-2})(-\pi i\tau)^{-s}\,ds.
$$
By the Mellin inversion formula again, this is equal to
$$
  \frac i{8\tau}\int_0^q\int_0^q\sum_{n=1}^\infty
  \left(\frac{-1}n\right)\frac{n^2q^{n/2}}{1-q^n}\frac{dq}q\frac{dq}q
 =\frac i{8\tau}\int_0^q\int_0^q
  \frac{\eta(\tau)^4\eta(4\tau)^8}{\eta(2\tau)^6}
  \frac{dq}q\frac{dq}q.
$$
This establishes \eqref{equation: cusp example 1}. In summary, what we
have shown is the following result.

\begin{Theorem} Let $\{a_n\}$ and $\{b_n\}$ be the
  sequences satisfying
$$
  (n+2)^2u_{n+2}-(32n^2+96n+76)u_{n+1}+(n+1)^2u_n=0 \quad
  (u_n=a_n\text{ or }b_n)
$$
with the initial values $a_0=1$, $a_1=12$, $b_0=0$, and $b_1=1$. Then
$$
  \left|\frac{b_n}{a_n}-\frac12L(2,\chi_{-1})\right|\ll\frac1{\log n}
$$
as $n\to\infty$.
\end{Theorem}

\noindent{\bf Example 4.} Let $f(\tau)=\eta(2\tau)^3\eta(6\tau)^3$ be
the cusp form of weight $3$ on $\Gamma_0(6)+\omega_3$ with character
$$
  \chi\begin{pmatrix}a&b\\6c&d\end{pmatrix}=(-1)^c\left(\frac
  d3\right), \qquad
  \chi\begin{pmatrix}3&-1\\6&-3\end{pmatrix}=i.
$$
We will find a rational approximation to $L(f,2)$.

Choose
\begin{equation} \label{equation: Example 4 t and A}
  t(\tau)=\frac{\eta(6\tau)^5\eta(2\tau)}{\eta(\tau)^5\eta(3\tau)},
  \qquad
  A(\tau)=\frac{\eta(\tau)^2\eta(3\tau)^3}{\eta(2\tau)\eta(6\tau)},
\end{equation}
where $t(\tau)$ is modular on the smaller group $\Gamma_0(6)$.
The differential equation satisfied by $t$ and $A$ is
$$
  (1+8t)^2(1+9t)\theta^2A+9t(1+8t)^2\theta A+2t(1+16t+72t^2)A=0.
$$
The singularities $0$, $-1/9$, $-1/8$, and $\infty$ correspond to the
cusps $\infty$, $1/2$, $1/3$, and $0$, respectively. Set
$g(t)=t/(1+8t)(1+9t)$ so that
$$
  \left(\frac{qdt/dq}t\right)^2\frac{g(t)}A=f(\tau).
$$
Let $A(t)=1+\sum_{n=1}^\infty a_nt^n$ and
$B(t)=t-7t^2+\sum_{n=3}^\infty b_nt^n$ be the solution of
$$
  (1+8t)^2(1+9t)\theta^2B+9t(1+8t)^2\theta B+2t(1+16t+72t^2)B=t+8t^2
$$
holomorphic at $t=0$. We now compute the limit of $b_n/a_n$.
For this purpose, we need to determine the behavior of
$$
  E(\tau)=\int_0^q\int_0^qf(\tau)\frac{dq}q\frac{dq}q
$$
near the cusp $1/2$. Choose
$$
  \sigma=\frac1{\sqrt 3}\begin{pmatrix}3&-2\\6&-3\end{pmatrix}
  \in\Gamma_0(6)+\omega_3,
$$
and consider $E(\sigma\tau)$.

Since $f(\sigma\tau)=i(2\sqrt 3\tau-\sqrt 3)^3f(\tau)$, according to
Lemma \ref{lemma: functional equation}, if $f(\tau)=\sum_{n=1}^\infty
c_nq^n$, then the function
$$
  L(s)=\sum_{n=1}^\infty(-1)^n\frac{c_n}{n^2}=-L(s,f)
$$
satisfies
$$
  \left(\frac{\pi}{\sqrt 3}\right)^{-s}\Gamma(s)L(s)
 =\left(\frac{\pi}{\sqrt 3}\right)^{s-3}\Gamma(3-s)L(3-s).
$$
Also, we have
$$
  E\left(-\frac1{2\sqrt 3\tau}+\frac12\right)
 =\frac1{2\pi i}\int_{3/2-\infty}^{3/2+i\infty}
  \Gamma(s)L(s+2)\left(\frac{\pi i}{\sqrt 3\tau}\right)^{-s}\,ds.
$$
The integrand has simple poles at $s=0$ and $s=-1$ with residues
$$
  \operatorname{Res}_{s=0}=L(2)=-L(2,f), \qquad
  \operatorname{Res}_{s=-1}=-L(1)\left(\frac{\pi i}{\sqrt3\tau}\right)
 =\frac i\tau L(2,f).
$$
Using the functional equation for $L(s)$ and the Mellin inversion
formula, we find
$$
  E\left(-\frac1{2\sqrt 3\tau}+\frac12\right)+L(2,f)
 =\frac i\tau\left(E\left(\frac\tau{2\sqrt 3}+\frac12\right)+L(2,f)\right).
$$
By Lemma \ref{lemma: alternative}, this amounts to
$$
  E(\sigma\tau)+L(2,f)=i(2\sqrt 3\tau-\sqrt3)^{-1}(E(\tau)+L(2,f)).
$$
Since $A(\tau)$ satisfies
$$
  A(\sigma\tau)=-i(2\sqrt 3\tau-\sqrt 3)A(\tau),
$$
we have
$$
  A(\sigma\tau)(E(\sigma\tau)+L(2,f))=A(\tau)(E(\tau)+L(2,f)).
$$
From this we see that the function $A(\tau)(E(\tau)+L(2,f))$, as a
function of $t$, is holomorphic at $t=-1/9$. Therefore, we have the
following approximation of $L(2,f)$.

\begin{Theorem} Let $f(\tau)=\eta(2\tau)^3\eta(6\tau)^3$ and $L(s,f)$
  be its associated $L$-function. Let $\{a_n\}$ and $\{b_n\}$ be the
  sequences satisfying
$$
  (n+3)^2u_{n+3}+(25n^2+109n+120)u_{n+2}+16(13n^2+35n+24)u_{n+1}
  +144(2n+1)^2u_n=0
$$
with the initial values
$$
  a_0=1, \quad a_1=-2, \quad a_2=10, \quad
  b_0=0, \quad b_1=1, \quad b_2=-7.
$$
Then we have
$$
  \left|\frac{b_n}{a_n}+L(2,f)\right|\ll(8/9)^n.
$$
\end{Theorem}

\noindent{\bf Remark.} Observe that the result in the above example
shows that the composition $t\mapsto\tau\mapsto
F(\tau)=A(\tau)(E(\tau)+L(2,f))$ is a single-valued analytic function
at $t=-1/9$. This in turn implies that the function $F(\tau)$ is
invariant under the substitution $\tau\mapsto\gamma\tau$ for any
$\gamma$ in the stabilizer subgroup of $1/2$ in $\Gamma_0(6)$. This
observation leads us to some interesting identities for the values of
$L$-functions. We will discuss these identities in the next section.
\medskip

\noindent{\bf Example 5.} Consider the differential equation
$$
  \theta^3A-8t(2\theta+1)(2\theta^2+2\theta+1)A+256t^2(\theta+1)^3A=0.
$$
This is equation ($\beta$) in \cite{ASZ}. We now determine its Ap\'ery
limit.

The differential equation is satisfied by
$$
  t(\tau)=\frac{\eta(\tau)^8\eta(4\tau)^{16}}{\eta(2\tau)^{24}},
  \qquad
  A(\tau)=\frac{\eta(2\tau)^{20}}{\eta(\tau)^8\eta(4\tau)^8}.
$$
They are modular on $\Gamma_0(4)$. The singularities $0,~1/16,~\infty$
correspond to the cusps $\infty,~0,~1/2$, respectively. The solution
of the inhomogeneous differential equation
$$
  \theta^3B-8t(2\theta+1)(2\theta^2+2\theta+1)B+256t^2(\theta+1)^3B=t
$$
is $B(\tau)=A(\tau)E(\tau)$, where
$$
  E(\tau)=\frac1{240}\int_0^q\int_0^q\int_0^q\left(
  E_4(\tau)-17E_4(2\tau)+16E_4(4\tau)\right)\frac{dq}q\frac{dq}q
  \frac{dq}q.
$$
We have
$$
  A(-1/4\tau)=-4\tau^2A(\tau).
$$
To determine the Ap\'ery limit, we should study how $E(\tau)$
transforms under the substitution $\tau\mapsto-1/4\tau$.

As usual, we apply the Mellin inversion formula and write
$$
  E(-1/4\tau)=\frac1{2\pi i}\int_{3/2-i\infty}^{3/2+i\infty}
  \Gamma(s)L(s+3)\left(\frac{2\pi i}{4\tau}\right)^{-s}\,ds,
$$
where
$$
  L(s)=\zeta(s)\zeta(s-3)(1-17\cdot 2^{-s}+16\cdot 2^{-2s})
$$
satisfies
$$
  \pi^{-s}\Gamma(s)L(s)=\pi^{s-4}\Gamma(4-s)L(4-s).
$$
Moving the line of integration to $\Re s=-7/2$, counting the residues,
making a change of variable $s\to -2-s$, using the functional equation
above, and invoking the Mellin inversion formula, we obtain
$$
  E(-1/4\tau)=\frac7{16}\zeta(3)-\frac{\pi^3i}{64\tau}
  -\frac7{64\tau^2}\zeta(3)+\frac1{4\tau^2}E(\tau).
$$
From this we deduce the following result.

\begin{Theorem} Let $\{a_n\}$ and $\{b_n\}$ be the sequences
  satisfying
$$
  (n+2)^3u_{n+2}-8(2n+3)(2n^2+6n+5)u_{n+1}+256n^3u_n=0 \quad
  (u_n=a_n\text{ or }b_n)
$$
with the initial values $a_0=1$, $a_1=8$, $b_0=0$, and $b_1=1$.
Then we have
$$
  \left|\frac{b_n}{a_n}-\frac7{16}\zeta(3)\right|\ll\frac1{\log n}
$$
as $n\to\infty$.
\end{Theorem}
\end{section}

\begin{section}{Identities involving $L$-values}
\label{section: L values}
Let the functions $t(\tau)$, $A(\tau)$, $f(\tau)$, and $E(\tau)$ be
given as in Example 4. Following the remark at the end of the example,
we know that the function $F(\tau)=A(\tau)(E(\tau)+L(2,f))$ satisfies
$F(\gamma^n\tau)=F(\tau)$ for all integers $n$, where
$\gamma=\left(\begin{smallmatrix}7&-3\\12&-5\end{smallmatrix}\right)$
is a generator of the stabilizer subgroup of the cusp $1/2$ in
$\Gamma_0(6)$. Now we have $A(\gamma\tau)=(12\tau-5)A(\tau)$. Thus,
the function $E(\tau)$ satisfies
$E(\gamma\tau)+L(2,f)=(12\tau-5)^{-1}(E(\tau)+L(2,f))$. By Lemma
\ref{lemma: alternative}, this is equivalent to
\begin{equation} \label{equation: remark 1}
  E\left(\frac\tau{12}+\frac7{12}\right)+L(2,f)=-\tau\left(
  E\left(-\frac1{12\tau}+\frac5{12}\right)+L(2,f)\right).
\end{equation}

On the other hand, we have $f(\gamma\tau)=(12\tau-5)^3f(\tau)$.
By Lemma \ref{lemma: functional equation}, letting $c_n$ denote the
Fourier coefficients of $f(\tau)=\sum_{n=1}^\infty c_nq^n$, the pair
of functions
$$
  L(s)=\sum_{n=1}^\infty\frac{c_n}{n^2}e^{14\pi i/12}, \qquad
  L^\ast(s)=\sum_{n=1}^\infty\frac{c_n}{n^2}e^{10\pi i/12}
$$
satisfy
$$
  \left(\frac\pi 6\right)^{-s}\Gamma(s)L(s)=i^3
  \left(\frac\pi 6\right)^{s-3}\Gamma(3-s)L^\ast(3-s).
$$
Then, arguing as before, we have
\begin{equation} \label{equation: remark 2}
  E\left(-\frac1{12\tau}+\frac5{12}\right)
  =L^\ast(2)-\frac1\tau L(2)-\frac1\tau
  E\left(\frac\tau{12}+\frac7{12}\right).
\end{equation}
Comparing \eqref{equation: remark 1} and \eqref{equation: remark 2},
we find
$$
  \sum_{n=1}^\infty\frac{c_n}{n^2}e^{14\pi in/12}
 =\sum_{n=1}^\infty\frac{c_n}{n^2}e^{10\pi in/12}
 =-\sum_{n=1}^\infty\frac{c_n}{n^2}.
$$

Now we observe that the deduction of \eqref{equation: remark 1} and
\eqref{equation: remark 2} is really independent of differential
equations and rational approximations. Thus, we expect that identities
of this type must exist in a general situation.

\begin{Lemma} \label{lemma: L values} Let $\Gamma=\Gamma_0(N)$ or
  $\Gamma_1(N)$. Let $\alpha=p/q\neq\infty$ be a cusp of $\Gamma$. Let
  $f(\tau)=\sum_{n=1}^\infty c_nq^n$ be a holomorphic modular form of
  weight $3$ on $\Gamma$ with character $\chi$. Assume that $f(\tau)$
  vanishes at $\infty$ and $\alpha$. Then for any
  $\gamma=\left(\begin{smallmatrix}a&b\\c&d\end{smallmatrix}\right)$
  with $c>0$ in the stabilizer subgroup of $\alpha$ in $\Gamma$ such
  that $\chi(\gamma)=1$, we have
  $$
    L(2)=L^\ast(2)=L_\alpha(2),
  $$
  where the $L$-functions are defined by
  $$
    L(s)=\sum_{n=1}^\infty\frac{c_n}{n^s}e^{2\pi ina/c}, \quad
    L^\ast(s)=\sum_{n=1}^\infty\frac{c_n}{n^s}e^{-2\pi ind/c}, \quad
    L_\alpha(s)=\sum_{n=1}^\infty\frac{c_n}{n^s}e^{2\pi in\alpha}
  $$
  for $\Re s>3$ and then continued analytically to the whole complex
  plane.
\end{Lemma}

\begin{proof} Let $\CC$ be the path consisting of the vertical line
  from $-d/c+i\infty$ to $-d/c$, the semi-circle from $-d/c$ to
  $\alpha$, and the vertical line from $\alpha$ to $\alpha+i\infty$.
  Consider the integral
  $$
  \int_\CC(c\tau+d)f(\tau)\,d\tau.
  $$
  Since $-d/c$ is equivalent to $\infty$ under $\Gamma$, by the
  assumption that $f(\tau)$ vanishes at $\infty$ and $\alpha$, the
  integral is convergent. Therefore, it is equal to $0$ because
  $f(\tau)$ is assumed to be a holomorphic modular form. Now we make a
  change of variable $\tau\to\gamma^{-1}\tau$ in the integral from
  $-d/c$ to $\alpha$. We have
  $$
    \int_{-d/c}^\alpha(c\tau+d)f(\tau)\,d\tau
   =\int_{i\infty}^\alpha(c\gamma^{-1}\tau+d)f(\gamma\tau)
    \frac{d\tau}{(-c\tau+a)^2}
   =\int_{i\infty}^\alpha f(\tau)\,d\tau.
  $$
  It follows that
  \begin{equation} \label{equation: theorem 1}
    \int_{-d/c}^{-d/c+i\infty}(c\tau+d)f(\tau)\,d\tau
   =\int_\alpha^{\alpha+i\infty}(c\tau+d-1)f(\tau)\,d\tau.
  \end{equation}
  Now recall that a stabilizer $\gamma$ of the cusp $\alpha=p/q$ takes
  the from
  $$
    \begin{pmatrix} 1+pqm & -p^2m\\ q^2m& 1-pqm\end{pmatrix}
  $$
  for some integer $m$. Hence, \eqref{equation: theorem 1} can be
  written as
  $$
    \int_{-d/c}^{-d/c+i\infty}(c\tau+d)f(\tau)\,d\tau
   =\int_{\alpha}^{\alpha+i\infty}(c\tau-c\alpha)f(\tau)\,d\tau.
  $$
  Finally, the two sides of the above identity are equal to
  $$
    \frac{ic}{4\pi^2}L^\ast(2), \qquad\frac{ic}{4\pi^2}L_\alpha(2),
  $$
  respectively. This gives $L^\ast(2)=L_\alpha(2)$. By a similar
  argument, we can also show that $L(2)=L_\alpha(2)$. This proves the
  lemma.
\end{proof}


\noindent{\bf Example 6.} Consider the differential equation
$$
  (1-11t-t^2)\theta^2A-(11t+2t^2)\theta A-(3t+t^2)A=0,
$$
appearing in Ap\'ery's sequence for $\zeta(2)$, where $t$ and $A$ are
given by \eqref{equation: zeta(2) t} and \eqref{equation: zeta(2) A}.
Using Lemma \ref{lemma: B(t)}, we find a solution for the
inhomogeneous differential equation
$$
  (1-11t-t^2)\theta^2B-(11t+2t^2)\theta B-(3t+t^2)B=t
$$
is given by $B(t)=A(t)E(t)$, where
\begin{equation} \label{equation: Example 6-1}
  E(t)=\int_0^q\int_0^q\left(\sum_{n\equiv 1\Mod 5}
  -2\sum_{n\equiv 2\Mod 5}+2\sum_{n\equiv 3\Mod 5}
  -\sum_{n\equiv 4\Mod 5}\right)\frac{n^2q^n}{1-q^n}\frac{dq}q
  \frac{dq}q.
\end{equation}
In \cite{Beukers}, Beukers mentioned that it can be shown that
$$
  A\left(\frac{\tau}{5\tau+1}\right)\left(
  E\left(\frac{\tau}{5\tau+1}\right)-\frac{\zeta(2)}5\right)
 =A(\tau)\left(E(\tau)-\frac{\zeta(2)}5\right),
$$
but did not give details. Here, we will prove this transformation
formula using Lemma \ref{lemma: L values}.

Since $A(\tau/(5\tau+1))=(5\tau+1)A(\tau)$, by Lemma \ref{lemma:
  alternative}, it suffices to prove that
\begin{equation} \label{equation: Example 6-2}
  E\left(-\frac1{5\tau}-\frac15\right)-\frac{\zeta(2)}5=
  -\frac1\tau\left(E\left(\frac\tau5+\frac15\right)-\frac{\zeta(2)}5
  \right).
\end{equation}
Let $f(\tau)=\sum_{n=1}c_nq^n$ be the modular form of weight $3$ in
the definition of $E(\tau)$ in \eqref{equation: Example 6-1}. It
satisfies $f(\tau/(5\tau+1))=(5\tau+1)^3f(\tau)$. Now, by the usual
argument, we find
$$
  E\left(-\frac1{5\tau}-\frac15\right)=L^\ast(2)+\frac{L(2)}\tau
 -\frac1\tau E\left(\frac\tau5+\frac15\right),
$$
where $L(s)=\sum c_ne^{2\pi in/5}/n^s$ and $L^\ast=\sum c_nn^{-2\pi
  in/5}/n^2$. Now applying Lemma \ref{lemma: L values} with
$\left(\begin{smallmatrix}1&0\\5&1\end{smallmatrix}\right)$ that fixes
$0$, we have
$$
  L(2)=L^\ast(2)=L(2,f)=\sum_{n=1}^\infty\frac{c_n}{n^2},
$$
which is shown to be $\zeta(2)/5$ in \cite{Beukers}. This gives the
desired equality \eqref{equation: Example 6-2}.
\medskip

\begin{Corollary} Let $\sum_{n=1}^\infty c_nq^n$ be the Fourier
  expansion of $f(\tau)=\eta(2\tau)^3\eta(6\tau)^3$. Let
  $L(s,f)=\sum_{n=1}^\infty c_n/n^s$ be the $L$-function associated
  with $f$. Then we have
  $$
    \sum_{n\equiv 1\Mod 12}\frac{c_n}{n^2}=\frac{2+\sqrt 3}3L(2,f),
    \qquad
    \sum_{n\equiv 7\Mod 12}\frac{c_n}{n^2}=\frac{2-\sqrt 3}3L(2,f).
  $$
\end{Corollary}

\begin{proof} Obviously, $c_n$ is $0$ whenever $n$ is even. Also,
  using Jacobi's triple product identity, we have
$$
  \eta(2\tau)^3\eta(6\tau)^3
 =\left(\sum_{n=0}^\infty(-1)^n(2n+1)q^{(2n+1)^2/4}\right)
  \left(\sum_{n=0}^\infty(-1)^n(2n+1)q^{3(2n+1)^2/4}\right).
$$
From this we see that $c_n=0$ whenever $n\equiv 5,11\Mod 12$.
Thus, letting $L_1=\sum_{n\equiv 1\Mod 12}c_n/n^2$, $L_2=\sum_{n\equiv
  7\Mod 12}c_n/n^2$, and $L_3=\sum_{3|n}c_n/n^2$, we have
$L(2,f)=L_1+L_2+L_3$. Moreover, since $c_{3n}=c_3c_n$ for all $n$, we
have $L_3=-L(2,f)/3$. It follows that
\begin{equation} \label{equation: corollary 1-1}
  L_1+L_2=\frac43L(2,f).
\end{equation}

We now apply Lemma \ref{lemma: L values} with the matrix
$\left(\begin{smallmatrix}1&0\\12&1\end{smallmatrix}\right)$, which
fixes the cusp $0$. Lemma \ref{lemma: L values} gives
$$
  e^{2\pi i/12}(L_1-L_2)+\sum_{n=1}^\infty\frac{c_{3n}}{(3n)^2}i^n
 =L(2,f).
$$
Comparing the real parts, we find
\begin{equation} \label{equation: corollary 1-2}
  L_1-L_2=\frac2{\sqrt 3}L(2,f).
\end{equation}
Combining \eqref{equation: corollary 1-1} and \eqref{equation:
  corollary 1-2}, we obtain the claimed identities.
\end{proof}

\begin{Corollary} Let $\sum_{n=1}^\infty c_nq^n$ be the Fourier
  expansion of $f(\tau)=\eta(4\tau)^6$. Let $L(s,f)=\sum_{n=1}^\infty
  c_n/n^s$ be the $L$-function associated with $f$. Then we have
  \begin{equation*}
  \begin{split}
    \sum_{n\equiv 1\Mod 16}\frac{c_n}{n^2}-\sum_{n\equiv 9\Mod 16}
    \frac{c_n}{n^2}&=\cos\frac{\pi}{8}L(2,f), \\
    \sum_{n\equiv 5\Mod 16}\frac{c_n}{n^2}-\sum_{n\equiv 13\Mod 16}
    \frac{c_n}{n^2}&=-\sin\frac{\pi}{8}L(2,f).
  \end{split}
  \end{equation*}
\end{Corollary}

\begin{proof} Clearly, the coefficients $c_n$ vanishes for all
  $n$ not congruent to $1$ modulo $4$. Let $L_i$, $i=0,\dots,3$
  denotes $\sum_{n\equiv 4i+1\Mod 16}^\infty c_n/n^2$. Apply Lemma
  \ref{lemma: L values} with the matrix
  $\left(\begin{smallmatrix}1&0\\16&1\end{smallmatrix}\right)$, which
  fixes the cusps $0$. We obtain
  $$
    e^{2\pi i/16}(L_0-L_2)+e^{10\pi i/16}(L_1-L_3)=L(2,f).
  $$
  Comparing the real parts and the imaginary parts of the two sides,
  we find
  $$
    \cos\frac{\pi}{8}(L_0-L_2)-\sin\frac{\pi}{8}(L_1-L_3)=L(2,f),
    \quad\sin\frac{\pi}{8}(L_0-L_2)+\cos\frac{\pi}{8}(L_1-L_3)=0.
  $$
  From this, we obtain the claimed identities.
\end{proof}
\end{section}

\end{document}